\documentclass[a4paper,12pt]{article}
\usepackage[utf8]{inputenc} 
\usepackage{amssymb,amsfonts}
\usepackage{amsfonts,amsmath,amssymb,amscd,latexsym}
\usepackage{srcltx}

\linethickness{0.4pt}
\ifx\plotpoint\undefined\newsavebox{\plotpoint}\fi 

\newcounter{ppp}
\newcommand{\bi}{\bibitem}
\newcommand{\nb}{\newblock}

\newcommand{\be}[1]{\begin{equation}\label{#1}}
\newcommand{\ee}{\end{equation}}

\newcommand{\la}{\langle\,}
\newcommand{\ra}{\,\rangle}
\newcommand{\ve}{\varepsilon}
\newcommand{\prf}{{\bf Proof.}\ }

\newcommand{\pp}{{\cal P}}

\newcommand{\dd}{{\cal D}}

\newcommand{\topp}{\mathop{\mbox{\bf top}}}
\newcommand{\bott}{\mathop{\mbox{\bf bot}}}

\newtheorem{thm}{\quad Theorem}
\newtheorem{lm}{\quad Lemma}

\newtheorem{prop}{\quad Proposition}

\title{On zero-measured subsets of Thompson's group $F$}
\author{\vspace{2ex}
	Victor Guba\thanks{This work is supported by the Russian Science Foundation, project no. 23-21-00289.}\\
	Vologda State University,\\
	15 Lenin Street,\\
	Vologda\\
	Russia\\
	160600\\
	E-mail: gubavs{@}vogu35.ru}
\date{}

\begin{document}
	
	\maketitle
	
	\begin{abstract}
		
A (discrete) group is called amenable whenever there exists a finitely additive right invariant probablity measure on it.
For Thompson's group $F$ the problem whether it is amenable is a long-standing open question. We consider presentation of $F$ in terms of non-spherical semigroup diagrams. There is a natural partition of $F$ into 7 parts in terms of these diagrams. We show that for any measure with the above properties on $F$, all but one of these sets have zero measure. This helps to clarify the structure of Folner sets in $F$ provided the group is amenable.
	
	\end{abstract}
	\vspace{5ex}
	
Here is the outline of the paper. Basic definitions are contained below in the text.
\vspace{1ex}

In Section~\ref{afz} we discuss definitions of amenable groups and introduce the concept of a zero-measured subset. We show that elements of such subsets can be excluded from Folner sets. In Section~\ref{thcd} we recall some basic facts about Thompson's group $F$ with respect to the amenability problem for it. We mention that the construction due to Belk and Brown does not give the optimal estimate for the density of finite subgraphs in it. Thus the assumption that $F$ may be amenable becomes more truthful. This makes actual the approach of constructing Folner sets for $F$. In this case it is useful to know what subsets in $F$ will be zero-measured.

We represent elements of $F$ by canonical semigroup diagrams accoring to their normal forms. According to this structure, we partition $F$ into 7 subsets. In Section~\ref{part7}, we show that all but one of these subsets are zero-measured. So they can be excluded from sequences of Folner sets in the process of their constructing.

\section{Amenability, Folner sets, and zero-measured sets}
\label{afz}

A (discrete) group $G$ is called {\em amenable\/} whenever there exists a finitely additive right invariant probability measure on $G$. Formally, this is a mapping $\mu\colon{\cal P}(G)\to[0,1]$ satisfying the following properties:

\begin{itemize}
\item 
 $\mu(A\cup B)=\mu(A)+\mu(B)$ for any disjoint subsets $A,B\subseteq G$,
\item
$\mu(Ag)=\mu(A)$ for any $A\subseteq G$, $g\in G$,
\item
$\mu(G)=1$.
\end{itemize}

Equivalent definitions of amenability come if we claim left or two-sided invariance of the measure, that is, the condition
$\mu(gA)=\mu(A)$ or $\mu(Ag)=\mu(gA)=\mu(A)$ (for all $A\subseteq G$, $g\in G)$. The proof can be found in \cite{GrL}. Here we will work with right invariant measures.
\vspace{1ex}

A subset $A$ of a group $G$ is called {\em zero-measured} whenever for any finitely additive right invariant probability measure $\mu$ on $G$ one has $\mu(A)=0$.
\vspace{1ex}

Trivial examples of zero-measured sets in infinite groups are finite subsets. Also we mention that if a group $G$ is non-amenable then no measures with the above properties exist on such a group. In this case we formally have that all subsets (including $G$ itself) are zero-measured.
\vspace{1ex}

The amenability property can be also stated in terms of Cayley graphs of groups. Here we restrict ourselves to the case when a group $G$ is finitely generated.
\vspace{1ex}

By the {\em density} of a finite graph we mean its average vertex degree. More precisely, let $v_1$, \dots, $v_k$ be all vertices of the graph $\Gamma$. Let $\deg_\Gamma(v)$ denote the degree of a vertex $v$ in $\Gamma$, that is, the number
of oriented edges of $\Gamma$ that come out of $v$. Then

\be{dgrform}
\delta(\Gamma)=\frac{\deg_\Gamma(v_1)+\cdots+\deg_\Gamma(v_k)}k
\ee
is called the density of $\Gamma$.

For the Cayley graph of a group $G$ with $m$ generators, it is known that $G$ is amenable if and only if the supremum of densities of its finite subgraphs has its maximum value $2m$. This is essentially a Folner criterion of amenability coming from~\cite{Fol}. This property does not depend on the choice of a finite generating set.
\vspace{1ex}

Let $A_n$ ($n\ge1$) be a sequence of finite nonemty subsets in $G$. We assume a finite $m$-generating set for $G$. Now these sets can be regarded as subgraphs in the Cayley graph of $G$ in these generators. Let $\delta$ denote the density of a finite graph. The family $\{A_n\mid n\ge1\}$ is called the {\em Folner sequence} whenever
\be{folseq}
\sup\limits_{n}\delta(A_n)=2m.
\ee

One can replace $\sup$ by $\lim\limits_{n\to\infty}$ if necessary. The members of the sequence are briefy called {\em Folner sets}.
\vspace{1ex}

Given a Folner sequence on $G$, it is easy to construct a finitely additive right invariant probablity measure $\mu$ on $G$ in the following way. Let $Z\subseteq G$ be any subset. Choose any non-principal ultrafilter ${\cal U}$ on $\mathbb N$. Then define $\mu(Z)$ as the limit of the quotient $\frac{|A_n\cap Z|}{|A_n|}$ along the ultrafilter, where $|\cdot|$ denotes the cardinality of a finte set:
\be{muz}
\mu(Z)=\lim\limits_{\cal U}\frac{|A_n\cap Z|}{|A_n|}.
\ee
It is easy to verify that $\mu$ will be a finitely additive right invariant probablity measure. Notice that the invariance has to be checked only for the generators of the group. According to the definition of Folner sets, they are ``almost invariant'' with respect to right shifts by generators. This means that for any group generator $x$, the cardinality of the symmetric difference $|A_nx\,\Delta\,A_n|$ is ``small'' with respect to $|A_n|$ as $n\to\infty$. This implies that the limit for $Zx$ in~(\ref{muz}) will remain the same.
\vspace{1ex}

Now let $A_n$ ($n\ge1$) be a Folner sequence and $A$ be a zero-measured subset in $G$. We are going to remove elements in $A$ from Folner sets. This gives us a new sequence $B_n=A_n\setminus A$. It can contain some empty sets. If we remove them from the sequence, then we get a new infinite family of nonempty finite subsets. We claim that the new sequence will be again Folner. We briefly state it as follows.

\begin{prop}
\label{remzero}
Removing elements of a zero-measured set from Folner sets, gives us a new sequence of Folner sets.
\end{prop}

First we need the following elementary

\begin{lm}
\label{del}
Let $\Gamma$ be a finite graph with density $\delta$ with $n$ vertices, where the degree of every vertex does not exceed $2m$. Suppose that we delete $k < \lambda n$ vertices from $\Gamma$ together with edges incident to them, where $\lambda < 1$. Then the density of the graph $\Gamma'$ obtained as a result is at least $\delta-2m\lambda$.
\end{lm}

$\prf$ Let $S$ be the sum of degrees of vertices in $\Gamma$. Then $\delta=\frac{S}n$ by definition. When we delete a vertex together with edges incident to it, we can loose at most $2m$ for the sum of degrees of the vertices in the rest. Hence the sum of degrees of vertices in $\Gamma'$ will be $S'\ge S-2mk$. For the density $\delta'$ of $\Gamma'$ we have $\delta'=\frac{S'}{n-k}\ge\frac{S'}n\ge\frac{S-2mk}n=\delta-2m\frac{k}n > \delta-2m\lambda$. The proof is complete.
\vspace{1ex}

{\bf Proof of Proposition~\ref{remzero}.} Let $A$ be zero-measured. First of all, we establish that  $\lim\limits_{n\to\infty}\frac{|A_n\cap A|}{|A_n|}=0$. Indeed, if this is not the case, then  $\limsup\limits_{n\to\infty}\frac{|A_n\cap A|}{|A_n|} > 0$. Then there exists an increasing infinite subsequence $n_1 < n_2 < \cdots < n_k < \cdots$ such that $\alpha=\lim\limits_{k\to\infty}\frac{|A_{n_k}\cap A|}{|A_{n_k}|} > 0$. Choose a non-principal ultrafilter ${\cal U}$ on $\mathbb N$ such that the set $\{n_k\mid k\in\mathbb N\}$ belongs to it. Let $\mu$ be a finitely additive right invariant probablity measure on $G$ constructed with respect to~(\ref{muz}). By definition, one has $\mu(A)=\lim\limits_{\cal U}\frac{|A_n\cap A|}{|A_n|}=\alpha > 0$. This is a contradiction.

Let $\varepsilon > 0$. Find a number $n_0$ such that $|A_n\cap A| < \frac{\varepsilon}{4m}|A_n|$ for $n\ge n_0$. Since $\sup\limits_{n}\delta(A_n)=2m$ for a Folner sequence, one can find a number $n > n_0$ such that $\delta(A_n) > 2m-\frac{\varepsilon}2$ (notice that the density never equals $2m$ since the group is infinite). Let $B_n=A_n\setminus A$. To get $B_n$ from $A_n$, we delete no more than $\lambda|A_n|$ vertices from $A_n$, where $\lambda=\frac{\varepsilon}{4m}$. According to Lemma~\ref{del}, we have $\delta(B_n) > \delta(A_n)-2m\lambda > 2m-\frac{\varepsilon}2-2m\frac{\varepsilon}{4m}=2m-\varepsilon$. Therefore, $\sup\limits_{n}\delta(B_n)=2m$ whence $B_n$ is a Folner sequence.

The proof is complete.
\vspace{1ex}

The r\^ole of Proposition~\ref{remzero} is the following: if we are going to construct a Folner sequence for a group $G$ and if we know that some set $A\subset G$ is zero-measured, then we can avoid elements in $A$ as members of the Folner sets.

\section{Thompson's group $F$ and canonical diagrams}
\label{thcd}

R.\,Thompson's group $F$ is given by the following infinite group presentation:

\be{xinf}
\la x_0,x_1,x_2,\ldots\mid x_j{x_i}=x_ix_{j+1}\ (i<j)\,\ra.
\ee

In many papers the same group is defined in terms of piecewise-linear functions. There are several ways to do that. All these definitions are equivalent. One of the basic references to the properties of this group is~\cite{CFP}. See also our recent survey~\cite{Gu22s}.

This group was discovered by Richard J. Thompson in the 60s. Some important properties of this group were obtained 
in~\cite{Bro,BG}. In~\cite{BS} it was proved that $F$ has no free non-abelian subgroups. It is easy to see that for any $n\ge2$, one has $x_n=x_0^{-(n-1)}x_1x_0^{n-1}$ so the group is generated by $x_0$, $x_1$. This generating set of the group is called {\em standard}. In these generators the group can be given by the following presentation with two defining relations:

\be{x0-1}
\la x_0,x_1\mid x_1^{x_0^2}=x_1^{x_0x_1},x_1^{x_0^3}=x_1^{x_0^2x_1}\ra,
\ee
where $a^b\leftrightharpoons b^{-1}ab$.

The famous problem about amenability of $F$ is still open. The question whether $F$ is amenable was asked by Ross Geoghegan in 1979; see~\cite{Geo,Ger87}. There is no common opinion on the answer: some of specialists in this area believe to non-amenability, some of them believe that the group is amenable. There is a number of papers with attempts to solve it in both directions. The author always believed in non-amenability of $F$ trying to prove this property. Now this belief is not so strong. Let us explain why.

A natural question is: what are the best known estimates to the density of finite subgraphs in the Cayley graph of $F$  
in standard generators $\{x_0,x_1\}$? In~\cite{Gu04} it was proved that the density approaches $3$; in the Addendum to the same paper it was shown that densities strictly exceed this value. An essential improvement was made by Belk and Brown~\cite{Be04,BB05}. They constructed a family of finite subgraphs whose densities approach $3{.}5$. There were many  attempts to improve this estimate. Several authors hypothesized that this construction was optimal, what would imply non-amenabilty of $F$. See~\cite{Bur16} and Conjecture 1 in our paper~\cite{Gu22}.

However, it turned out that this conjecture was false. Recently we proved in~\cite{Gu23} that there exist finite subgraphs in the Cayley graph of $F$ in generators $x_0$, $x_1$ with density strictly exceeding $3{.}5$. This makes amenability of $F$ more truthful. So in this paper we discuss some properties of possible Folner sets in $F$.

Notice that if $F$ is amenable and Folner sets exist for it, then it is known that they have a huge size. In~\cite{Moore13} it is shown that they grow as a tower of exponents.
\vspace{1ex}

Recall that any element in $F$ has the unique {\em normal form} in the infinite generating set. Namely, this is the following expression:

\be{nf}
x_{i_1}x_{i_2}\cdots x_{i_s}x_{j_t}^{-1}\cdots x_{j_2}^{-1}x_{j_1}^{-1},
\ee
where $s,t\ge0$, $0\le i_1\le i_2\le\cdots\le i_s$, $0\le j_1\le j_2 \le\cdots\le j_t$ and the following is true: if (\ref{nf}) contains both $x_i$ and $x_i^{-1}$ for some $i\ge0$, then it also contains $x_{i+1}$ or $x_{i+1}^{-1}$ (in particular, $i_s\ne j_t$).
\vspace{1ex}

It is also known that $F$ is a diagram group over the simplest semigroup presentaion $\la x\mid x^2=x\ra$. A detailed information on diagram groups can be found in~\cite{GbS}. Here we need to describe a modified version of this idea based on the representation of $F$ by non-spherical diagrams. More detailed explanation is conained in~\cite[Section 3]{Gu04}.
\vspace{1ex}

First we recall the concept of a semigroup diagram and introduce some notation. Let us consider the following example.
Let $\pp=\la a,b\mid aba=b,bab=a\ra$ be the semigroup presentation. It is easy to see by the following algebraic calculation
$$
a^5=a(bab)a(bab)a=(aba)(bab)(aba)=bab=a
$$
that the words $a^5$ and $a$ are equal modulo $\pp$. The same can be seen from the following picture

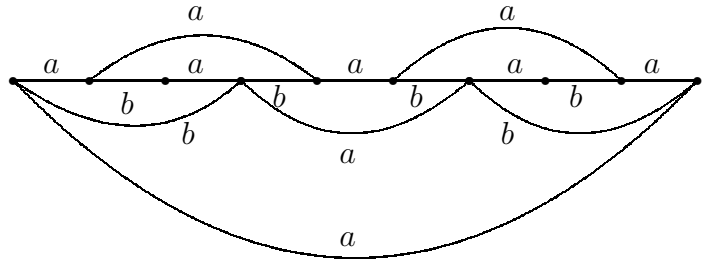
\begin{center}
	\begin{picture}(90.00,37.00)
		\put(00.00,23.00){\circle*{1.00}}
		\put(10.00,23.00){\circle*{1.00}}
		\put(20.00,23.00){\circle*{1.00}}
		\put(30.00,23.00){\circle*{1.00}}
		\put(30.00,23.00){\circle*{1.00}}
		\put(40.00,23.00){\circle*{1.00}}
		\put(50.00,23.00){\circle*{1.00}}
		\put(60.00,23.00){\circle*{1.00}}
		\put(60.00,23.00){\circle*{1.00}}
		\put(70.00,23.00){\circle*{1.00}}
		\put(80.00,23.00){\circle*{1.00}}
		\put(90.00,23.00){\circle*{1.00}}
		\put(00.00,23.00){\line(1,0){90.00}}
		\bezier{152}(10.00,23.00)(25.00,35.00)(40.00,23.00)
		\bezier{240}(50.00,23.00)(80.00,23.00)(50.00,23.00)
		\bezier{164}(50.00,23.00)(65.00,37.00)(80.00,23.00)
		\bezier{240}(0.00,23.00)(30.00,23.00)(0.00,23.00)
		\bezier{156}(0.00,23.00)(17.00,11.00)(30.00,23.00)
		\bezier{164}(30.00,23.00)(44.00,9.00)(60.00,23.00)
		\bezier{164}(60.00,23.00)(74.00,9.00)(90.00,23.00)
		\put(5.00,25.00){\makebox(0,0)[cc]{$a$}}
		\put(24.00,32.00){\makebox(0,0)[cc]{$a$}}
		\put(45.00,25.00){\makebox(0,0)[cc]{$a$}}
		\put(65.00,32.00){\makebox(0,0)[cc]{$a$}}
		\put(84.00,25.00){\makebox(0,0)[cc]{$a$}}
		\put(23.00,16.00){\makebox(0,0)[cc]{$b$}}
		\put(44.00,13.00){\makebox(0,0)[cc]{$a$}}
		\put(65.00,16.00){\makebox(0,0)[cc]{$b$}}
		\put(15.00,20.00){\makebox(0,0)[cc]{$b$}}
		\put(24.00,25.00){\makebox(0,0)[cc]{$a$}}
		\put(35.00,21.00){\makebox(0,0)[cc]{$b$}}
		\put(53.00,21.00){\makebox(0,0)[cc]{$b$}}
		\put(66.00,25.00){\makebox(0,0)[cc]{$a$}}
		\put(74.00,21.00){\makebox(0,0)[cc]{$b$}}
		\bezier{520}(0.00,23.00)(45.00,-24.00)(90.00,23.00)
		\put(44.00,2.00){\makebox(0,0)[cc]{$a$}}
	\end{picture}
\end{center}

This object is called a {\em diagram\/} $\Delta$ over the semigroup presentation $\pp$. It is a plane graph with $10$ vertices, $15$ (geometric) edges and $6$ cells. Each cell corresponds to an elementary transformation of a
word, that is, a transformation of the form $p\cdot u\cdot q\to p\cdot v\cdot q$, where $p$, $q$ are words (possibly,
empty), $u=v$ or $v=u$ belongs to the set of defining relations. The diagram $\Delta$ has the leftmost vertex denoted by $\iota(\Delta)$ and the rightmost vertex denoted by $\tau(\Delta)$. It also has the {\em top path\/} $\topp(\Delta)$ and the {\em bottom path\/} $\bott(\Delta)$ from $\iota(\Delta)$ to $\tau(\Delta)$. Each cell $\pi$ of a diagram can be
regarded as a diagram itself. The above functions $\iota$, $\tau$, $\topp$, $\bott$ can be applied to $\pi$ as well. We do not distinguish isotopic diagrams.

We say that $\Delta$ is a $(w_1,w_2)$-diagram whenever the label of its top path is $w_1$ and the label of its bottom path is $w_2$. In our example, we deal with an $(a^5,a)$-diagram. If we have two diagrams such that the bottom path of the first of them has the same label as the top path of the second, then we can naturally {\em concatenate\/} these diagrams by
identifying the bottom path of the first diagram with the top path of the second diagram. The result of the concatenation of a $(w_1,w_2)$-diagram and a $(w_2,w_3)$-diagram obviously is a $(w_1,w_3)$-diagram. We use the sign $\circ$ for the operation of concatenation. For any diagram $\Delta$ over $\pp$ one can consider its {\em mirror image\/} $\Delta^{-1}$ with respect to a horizontal axis. A diagram may have {\em dipoles\/}, that is, subdiagrams of the form $\pi\circ\pi^{-1}$,
where $\pi$ is a single cell. To {\em cancel\/} (or {\em reduce\/}) the dipole means to remove the common boundary of $\pi$ and $\pi^{-1}$ and then to identify $\topp(\pi)$ with $\bott(\pi^{-1})$. In any diagram, we can cancel all its dipoles, step by step. The result does not depend on the order of cancellations. A diagram is {\em reduced\/} whenever it has no dipoles. The operation of cancelling dipoles has an inverse operation called the {\em insertion\/} of a dipole. These operations induce an equivalence relation on the set of diagrams (two diagrams are {\em equivalent\/} whenever
one can go from one of them to the other by a finite sequence of cancelling/inserting dipoles). Each equivalence class contains exactly one reduced diagram.

For any nonempty word $w$, the set of all $(w,w)$-diagrams forms a monoid with the identity element $\ve(w)$ (the diagram with no cells). The operation $\circ$ naturally induces some operation on the set of equivalence classes of diagrams. This operation is called a {\em product\/} and equivalent diagrams are called {\em equal\/}. (The sign $\equiv$ will be
used to denote that two diagrams are isotopic.) So the set of all equivalence classes of $(w,w)$-diagrams forms a group that is called the {\em diagram group\/} over $\pp$ with {\em base\/} $w$. We denote this group by ${\cal D}(\pp,w)$. We can think of this group as of the set of all reduced $(w,w)$-diagrams. The group operation is the concatenation with cancelling all dipoles in the result. An inverse element of a diagram is its mirror image. We also need one more natural
operation on the set of diagrams. By the {\em sum\/} of two diagrams we mean the diagram obtained by identifying the rightmost vertex of the first summand with the leftmost vertex of the second summand. This operation is also associative. The sum of diagrams $\Delta_1$, $\Delta_2$ is denoted by $\Delta_1+\Delta_2$.

We have already noticed that the group $F$ is the diagram group over the simplest semigroup presentation $\pp=\la x\mid x^2=x\ra$ with base $x$ (as well as for any base $x^k$, where $k\ge1$). All these diagrams are {\em spherical\/}, that is, they are $(w,w)$-diagrams for some word $w$. Now we describe the following modification. 

Let $\Delta$ be any diagram over $\pp=\la x\mid x^2=x\ra$, not necessarily spherical. Let us add an infinite sequence of edges on the right of $\Delta$, each edge is labelled by $x$. This object will be called an {\em infinite diagram\/} over $\pp$. Note that it has finitely many cells. An infinite diagram that corresponds to $\Delta$ will be denoted by $\hat\Delta$. It has the leftmost vertex $\iota(\hat\Delta)$ and two distinguished infinite paths starting at
$\iota(\hat\Delta)$, both labelled by the infinite power of $x$. These paths will be denoted by $\topp(\hat\Delta)$ and $\bott(\hat\Delta)$, respectively. The concept of a dipole in an infinite diagram is defined as above. The
same concerns the operations of deleting/inserting a dipole, the equivalence relation induced by these operations, and so on. Any two infinite diagrams can be naturally concatenated (the bottom path of the first factor is identified with the top path of the second factor). We use the same sign $\circ$ for this concatenation. The operation $\circ$
gives the set of all infinite diagrams a monoid structure. The identity of it is the infinite diagram without cells denoted by $\ve$. As in the case of spherical diagrams, the operation of concatenation induces a group operation on the set of all equivalence classes of infinite diagrams. Thus we have a group. We shall denote it by $\hat\dd(\pp,x)$. (This makes
sense for any semigroup presentation $\pp$ in an alphabet of one letter. Notice that we can forget about the labels working with a one-letter alphabet.) It is easy to see that the group we have will be isomorphic to $F$. Indeed, let $X_i$ be the infinite diagram

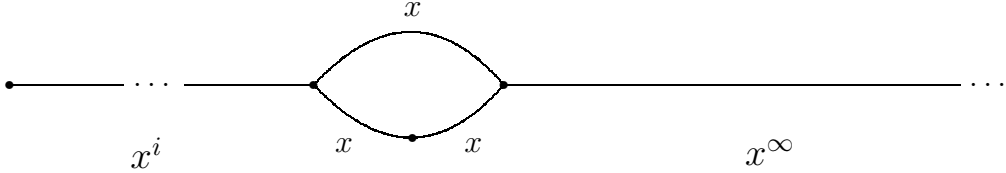
\begin{center}
	\begin{picture}(131.00,25.00)
		\put(2.00,11.00){\circle*{1.00}}
		\put(2.00,11.00){\line(1,0){15.00}}
		\put(21.00,11.00){\makebox(0,0)[cc]{$\dots$}}
		\put(25.00,11.00){\line(1,0){17.00}}
		\put(42.00,11.00){\circle*{1.00}}
		\put(67.00,11.00){\circle*{1.00}}
		\bezier{152}(42.00,11.00)(55.00,25.00)(67.00,11.00)
		\bezier{152}(42.00,11.00)(55.00,-3.00)(67.00,11.00)
		\put(67.00,11.00){\line(1,0){60.00}}
		\put(55.00,4.00){\circle*{1.00}}
		\put(131.00,11.00){\makebox(0,0)[cc]{$\dots$}}
		\put(55.00,21.00){\makebox(0,0)[cc]{$x$}}
		\put(46.00,3.00){\makebox(0,0)[cc]{$x$}}
		\put(63.00,3.00){\makebox(0,0)[cc]{$x$}}
		\put(20.00,2.00){\makebox(0,0)[cc]{\large$x^i$}}
		\put(102.00,2.00){\makebox(0,0)[cc]{\large$x^\infty$}}
	\end{picture}
\end{center}

\noindent
that consists of an $(x,x^2)$-cell, the finite path labelled by $x^i$ on the left of it and the infinite path labelled by the infinite power of $x$, on the right. By $X_i^{-1}$ we mean the mirror image of $X_i$ under the horizontal axis symmetry. Infinite diagrams of the form $X_i^{\pm1}$ ($i\ge0$) are called {\em atomic\/}. For any integers $j>i\ge0$, the diagram

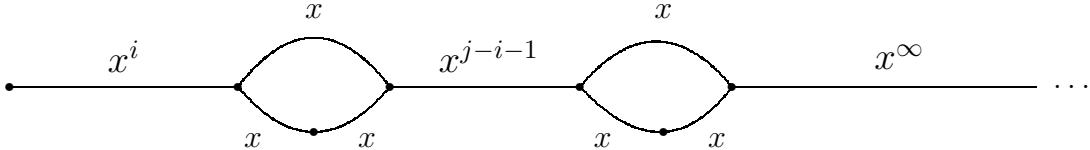
\begin{center}
	\begin{picture}(142.00,22.00)
		\put(2.00,9.00){\circle*{1.00}}
		\put(32.00,9.00){\circle*{1.00}}
		\put(52.00,9.00){\circle*{1.00}}
		\put(77.00,9.00){\circle*{1.00}}
		\put(97.00,9.00){\circle*{1.00}}
		\put(2.00,9.00){\line(1,0){30.00}}
		\put(52.00,9.00){\line(1,0){25.00}}
		\put(97.00,9.00){\line(1,0){40.00}}
		\bezier{132}(32.00,9.00)(42.00,22.00)(52.00,9.00)
		\bezier{132}(32.00,9.00)(42.00,-3.00)(52.00,9.00)
		\bezier{124}(77.00,9.00)(87.00,21.00)(97.00,9.00)
		\bezier{124}(77.00,9.00)(88.00,-3.00)(97.00,9.00)
		\put(42.00,3.00){\circle*{1.00}}
		\put(88.00,3.00){\circle*{1.00}}
		\put(42.00,19.00){\makebox(0,0)[cc]{$x$}}
		\put(88.00,19.00){\makebox(0,0)[cc]{$x$}}
		\put(34.00,2.00){\makebox(0,0)[cc]{$x$}}
		\put(49.00,2.00){\makebox(0,0)[cc]{$x$}}
		\put(80.00,2.00){\makebox(0,0)[cc]{$x$}}
		\put(95.00,2.00){\makebox(0,0)[cc]{$x$}}
		\put(17.00,13.00){\makebox(0,0)[cc]{\large$x^i$}}
		\put(65.00,13.00){\makebox(0,0)[cc]{\large$x^{j-i-1}$}}
		\put(119.00,13.00){\makebox(0,0)[cc]{\large$x^\infty$}}
		\put(142.00,9.00){\makebox(0,0)[cc]{$\dots$}}
	\end{picture}
\end{center}

\noindent
equals both $X_j\circ X_i$ and $X_i\circ X_{j+1}$. This means that we have a homomorphism from $F$ to the group of infinite diagrams. This homomorphism is onto because any infinite diagram is a concatenation of atomic diagrams. The homomorphism must be injective because all proper homomorphic images of $F$ are abelian \cite{CFP}. However, the group of infinite diagrams is not abelian since $X_1X_0=X_0X_2\ne X_0X_1$.

When we work with infinite diagrams, it is convenient to eliminate the infinite ``tail" on the right of each infinite diagram. An ordinary diagram over $\pp$ is called {\em canonical\/} whenever it has no dipoles and it is not a sum of a diagram and an edge. It is obvious that there is a one-to-one correspondence between the set of infinite diagrams without
dipoles and the set of canonical diagrams. So we may assume that each element of $F$ has a unique canonical representative. We have a group structure on the set of all canonical diagrams over $\pp$. Given an $(x^p,x^q)$-diagram $\Delta_1$ and an $(x^s,x^t)$-diagram $\Delta_2$, we multiply them as follows. If $q=s$, then we concatenate them. If $q<s$, then we concatenate $\Delta_1+\ve(x^{s-q})$ and $\Delta_2$. If $q>s$, then we concatenate $\Delta_1$ and $\Delta_2+\ve(x^{q-s})$.
After the concatenation, we reduce all dipoles in the result. Then we need to make the diagram canonical. This means that we have to delete the rightmost common suffix of the top and the bottom path of the diagram we have. The only exception is made for the diagram $\ve(x)$, the identity element of the group. This diagram is already canonical so we leave it as it is.

Given a normal form, it is easy to draw the corresponding diagram, and vice versa. The following example illustrates the diagram that corresponds to the element $g=x_0^3x_1x_3x_8x_{11}^2x_{12}x_{16}x_{17}x_{18}x_{17}^{-2}x_{11}^{-1}x_5^{-3}x_0^{-1}$ represented by its normal form:

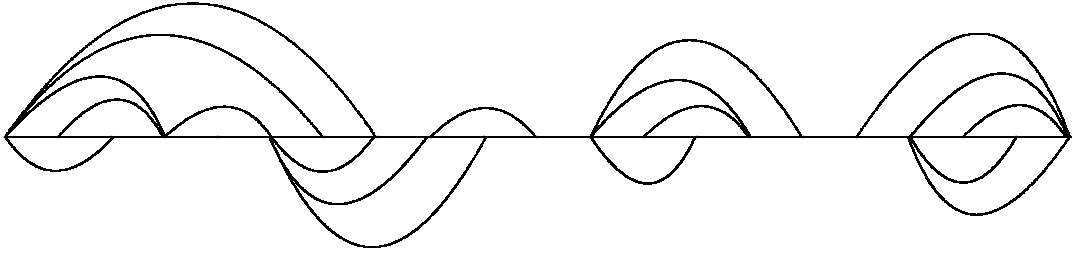
\begin{center}
\unitlength 1mm 
\linethickness{0.4pt}
\ifx\plotpoint\undefined\newsavebox{\plotpoint}\fi 
\begin{picture}(170.375,60.375)(0,0)
	\put(2.5,25){\line(1,0){140.2}}
	\put(2.5,25){\circle*{.5}}
	\put(9.5,25){\circle*{.5}}
	\put(16.5,25){\circle*{.5}}
	\put(23.5,25){\circle*{.5}}
	\put(30.5,25){\circle*{.5}}
	\put(37.5,25){\circle*{.5}}
	\put(44.5,25){\circle*{.5}}
	\put(51.5,25){\circle*{.5}}
	\put(58.5,25){\circle*{.5}}
	\put(65.5,25){\circle*{.5}}
	\put(72.5,25){\circle*{.5}}
	\put(79.5,25){\circle*{.5}}
	\put(86.5,25){\circle*{.5}}
	\put(93.5,25){\circle*{.5}}
	\put(100.5,25){\circle*{.5}}
	\put(107.5,25){\circle*{.5}}
	\put(114.5,25){\circle*{.5}}
	\put(121.5,25){\circle*{.5}}
	\put(128.5,25){\circle*{.5}}
	\put(135.5,25){\circle*{.5}}
	\put(142.25,25){\circle*{.5}}
	\qbezier(9.5,25)(18,34.875)(23.25,25)
	\qbezier(23.5,25)(32.125,33)(37.5,25)
	\qbezier(58.5,25)(66,32.625)(72.25,25)
	\qbezier(86.5,25)(94.75,33.125)(100.25,25)
	\qbezier(128.5,25)(136.5,33.375)(142,25)
	\qbezier(2.5,25)(16.875,41)(23.5,25)
	\qbezier(2.5,25)(22.5,52)(44.25,25)
	\qbezier(2.5,25)(27.375,60.375)(51.25,25)
	\qbezier(79.5,25)(92,40)(100.5,25)
	\qbezier(79.5,25)(91.625,50.625)(107.25,25)
	\qbezier(121.5,25)(135.125,41.75)(142.25,25)
	\qbezier(114.5,25)(132.625,52.375)(142.5,25)
	\qbezier(2.5,25)(8.625,16)(16.75,25)
	\qbezier(37.5,25)(44.25,15.75)(51,25)
	\qbezier(37.25,25)(44.875,7.125)(58,25)
	\qbezier(37.25,25)(50,-4.25)(65.75,25)
	\qbezier(79.5,25)(87.625,12.625)(93.25,25)
	\qbezier(121.5,25)(128.5,12.875)(135.5,25)
	\qbezier(121.25,25)(128.625,4.375)(142.5,25)
\end{picture}

\end{center}

\section{Partition of $F$ into 7 sets}
\label{part7}

We describe a partition of $F$ into 7 sets according to the structure of canonical diagrams representing these elements. 

Let $g$ be an element of $F$ given by its canonical diagram $\Delta$. For $i\ge0$ and $\epsilon=\pm1$ we say that $\Delta$ is {\em right divisible} by $X_i^{\epsilon}$ whenever $\Delta$ is a concatenation of some diagram $\Delta'$ and $X_i^{\epsilon}$. Equivalently, one can say that the diagram $\Delta\circ X_i^{-\epsilon}$ has a dipole. In this case the canonical diagram representing $gx_i^{-\epsilon}$ will be $\Delta'$.

By $D(\Delta)$ we denote the set of right divisors of $\Delta$ among the set $\{X_0,X_0^{-1},X_1,X_1^{-1}\}$. There are 7 cases for that. By ${\cal M}_i$ we denote the set of elements in $F$ whose canonical diagrams satisfy condition of the $i$th case ($1\le i\le7$). For an element $g\in F$ we will always denote by $\Delta=\Delta(g)$ its canonical diagram.
\vspace{1ex}

1) ${\cal M}_1=\{g\in F\mid D(\Delta)=\emptyset\}$. Examples of elements in this set: $e$, $x_2$, $x_1x_2^{-1}$.
\vspace{1ex}

2) ${\cal M}_2=\{g\in F\mid D(\Delta)=\{X_0^{-1}\}\}$. Examples of elements in this set: $x_0^{-1}$, $x_0x_1x_0^{-1}$, $x_2x_1^{-1}x_0^{-1}$.
\vspace{1ex}

3) ${\cal M}_3=\{g\in F\mid D(\Delta)=\{X_0\}\}$. Examples of elements in this set: $x_0$, $x_0x_3^{-1}$, $x_0^3x_2^{-1}$. Notice that if $\Delta$ is right divisible by $X_0$ then no other right divisors from the set $\{X_0^{-1},X_1,X_1^{-1}\}$ can occur. 
\vspace{1ex}

4) ${\cal M}_4=\{g\in F\mid D(\Delta)=\{X_1^{-1}\}\}$. Examples of elements in this set: $x_1^{-1}$, $x_0x_1^{-2}$, $x_3x_4x_1^{-1}$.
\vspace{1ex}

5) ${\cal M}_5=\{g\in F\mid D(\Delta)=\{X_1\}\}$. Examples of elements in this set: $x_1$, $x_0x_1$, $x_1^2x_4^{-1}$.
\vspace{1ex}

6) ${\cal M}_6=\{g\in F\mid D(\Delta)=\{X_0^{-1},X_1^{-1}\}\}$. Examples of elements in this set: $x_2^{-1}x_0^{-1}$, $x_1x_4^{-1}x_0^{-3}$, $x_1x_5x_3^{-2}x_0^{-2}$.
\vspace{1ex}

7) ${\cal M}_7=\{g\in F\mid D(\Delta)=\{X_0^{-1},X_1\}\}$. Examples of elements in this set: $x_2x_0^{-1}$, $x_0x_1x_2x_0^{-1}$, $x_1x_3x_0^{-2}$.
\vspace{1ex}

It is not hard to see that any canonical diagram satisfies exacly one of the 7 conditions, so we have a partition of $F$ into 7 subsets. Our aim is to show that if $\mu$ is a  finitely additive right invariant probability measure on $F$, then it is competely concentrated on ${\cal M}_6$. That is, all the other six sets are zero-measured.

To make the situation more visual, we add a picture showing typical cell structure of canonical diagrams $\Delta$ for these subsets. A cell that corresponds to $X_0^{\pm1}$, $X_1^{\pm1}$ is shown if and only if it belongs to $D(\Delta)$. 

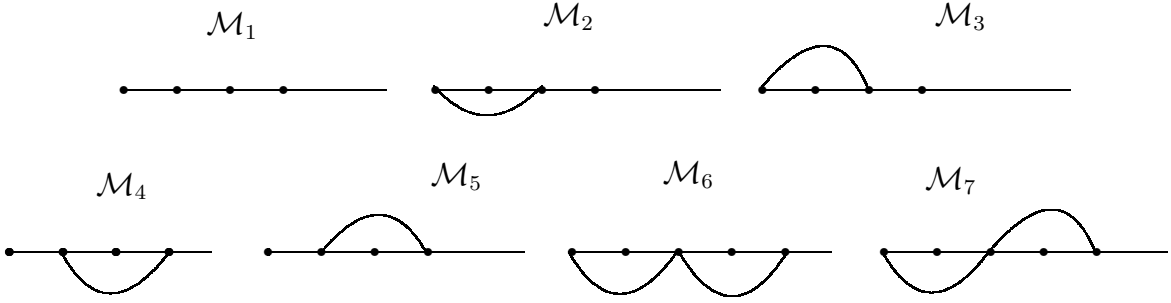
\begin{center}
\unitlength 1mm 
\linethickness{0.4pt}
\ifx\plotpoint\undefined\newsavebox{\plotpoint}\fi 

\begin{picture}(156.25,44.125)(0,0)
	\put(4,11.5){\circle*{1}}
	\put(4,11.5){\circle*{1}}
	\put(4,11.5){\circle*{1}}
	\put(4,11.5){\circle*{1}}
	\put(11,11.5){\circle*{1}}
	\put(11,11.5){\circle*{1}}
	\put(11,11.5){\circle*{1}}
	\put(11,11.5){\circle*{1}}
	\put(18,11.5){\circle*{1}}
	\put(18,11.5){\circle*{1}}
	\put(18,11.5){\circle*{1}}
	\put(18,11.5){\circle*{1}}
	\put(25,11.5){\circle*{1}}
	\put(25,11.5){\circle*{1}}
	\put(25,11.5){\circle*{1}}
	\put(25,11.5){\circle*{1}}
	\put(4,11.5){\line(1,0){26.5}}
	\put(4,11.5){\line(1,0){26.5}}
	\put(4,11.5){\line(1,0){26.5}}
	\put(4,11.5){\line(1,0){26.5}}
	\qbezier(10.75,11.5)(16.625,.5)(25,11.5)
	\qbezier(10.75,11.5)(16.625,.5)(25,11.5)
	\qbezier(10.75,11.5)(16.625,.5)(25,11.5)
	\qbezier(10.75,11.5)(16.625,.5)(25,11.5)
	\put(38,11.5){\circle*{1}}
	\put(38,12){\line(0,-1){1}}
	\put(38,11.5){\line(1,0){33.75}}
	\put(45,11.5){\circle*{1}}
	\put(52,11.5){\circle*{1}}
	\put(59,11.5){\circle*{1}}
	\qbezier(45,11.5)(53.125,21.375)(58.75,11.5)
	\put(78,11.5){\line(1,0){34}}
	\put(119,11.5){\line(1,0){37.25}}
	\put(78,11.5){\circle*{1}}
	\put(85,11.5){\circle*{1}}
	\put(92,11.5){\circle*{1}}
	\put(99,11.5){\circle*{1}}
	\put(106,11.5){\circle*{1}}
	\qbezier(77.5,11.5)(83.875,.375)(91.75,11.5)
	\qbezier(92,11.5)(98.875,-.25)(106.25,11.5)
	\put(119,11.5){\circle*{1}}
	\put(126,11.5){\circle*{1}}
	\put(133,11.5){\circle*{1}}
	\put(140,11.5){\circle*{1}}
	\put(147,11.5){\circle*{1}}
	\qbezier(118.75,11.5)(124.75,.75)(132.75,11.5)
	\qbezier(132.75,11.5)(142.25,22.75)(146.75,11.5)
	\put(18.75,20.25){\makebox(0,0)[cc]{${\cal M}_4$}}
	\put(62.75,21.5){\makebox(0,0)[cc]{${\cal M}_5$}}
	\put(93.25,21.75){\makebox(0,0)[cc]{${\cal M}_6$}}
	\put(127.75,21){\makebox(0,0)[cc]{${\cal M}_7$}}
	\put(19,33){\line(1,0){34.5}}
	\put(60,33){\line(1,0){37.5}}
	\put(103,33){\line(1,0){40.5}}
	\put(19,33){\circle*{1}}
	\put(26,33){\circle*{1}}
	\put(33,33){\circle*{1}}
	\put(40,33){\circle*{1}}
	\put(60,33){\circle*{1}}
	\put(67,33){\circle*{1}}
	\put(74,33){\circle*{1}}
	\put(81,33){\circle*{1}}
	\put(103,33){\circle*{1}}
	\put(110,33){\circle*{1}}
	\put(117,33){\circle*{1}}
	\put(124,33){\circle*{1}}
	\qbezier(59.75,33.5)(66.625,25.75)(74,33.5)
	\qbezier(103,33.5)(112.125,44.125)(116.75,33.5)
	\put(33.25,41.5){\makebox(0,0)[cc]{${\cal M}_1$}}
	\put(77.5,42.75){\makebox(0,0)[cc]{${\cal M}_2$}}
	\put(129,42.75){\makebox(0,0)[cc]{${\cal M}_3$}}
\end{picture}
\end{center}

The following Lemma shows what happens if we act on sets of the form ${\cal M}_i$ by right multiplication on some group generators $x_0^{\pm1}$, $x_1^{\pm1}$.

\begin{lm}
\label{rightm}

a) $({\cal M}_1\cup{\cal M}_3\cup{\cal M}_4\cup{\cal M}_5)x_0\subseteq{\cal M}_3$; 
b) $({\cal M}_2\cup{\cal M }_7)x_1\subseteq{\cal M}_7$;  c) ${\cal M}_7x_0^{-1}\subseteq{\cal M}_2$; ${\cal M}_3x_1^{-1}\subseteq{\cal M}_4$.
\end{lm}

\prf If $\Delta$ is not right divisible by $X_0^{-1}$ then $\Delta\circ X_0$ has no dipoles and it is right divisible by $X_0$. In this case $X_0$ is the only right divisor of the diagram, and the result belongs to ${\cal M}_3$. This implies a).

To prove b), let us take a diagram $\Delta$ from ${\cal M}_2$ or ${\cal M}_7$. Multiplying it on the right by $X_1$ means adding a cell of the form $x=x^2$ attached by its top path to the second bottom edge of $\Delta$. No dipoles occur in this case since if $\Delta$ belongs to ${\cal M}_2$ then it is not right divisible by $X_1^{-1}$. Therefore, we get to a canonical diagram from ${\cal M}_7$.

To establish c), we take a diagram $\Delta$ from ${\cal M}_7$. Multiplying it on the right by $X_0^{-1}$ means adding an arc connecting the first and the third vertex on the bottom. Clearly, no diploes occur after this operation, and the result belongs to ${\cal M}_2$.

As for d), if we take a diagram $\Delta$ from ${\cal M}_3$ and multiply it on the right by $X_1^{-1}$, then we add an arc connecting the 2nd and the 4th vertex on the bottom. Obviously, we get to ${\cal M}_4$.

The proof is complete.
\vspace{1ex}

Now our main result follows almost immediately.

\begin{thm}
\label{m6}
If $\mu$ is a finitely additive right invariant probability measure on $F$, then the subsets ${\cal M}_i$ $(1\le i\le 5)$ and ${\cal M}_7$ are zero-measured. That is, all measure (provided it exists) is concentrated on ${\cal M}_6$.
\end{thm}

\prf Part a) implies that $\mu({\cal M}_1)+\mu({\cal M}_3)+\mu({\cal M}_4)+\mu({\cal M}_5)=\mu({\cal M}_1\cup{\cal M}_3\cup{\cal M}_4\cup{\cal M}_5)=\mu(({\cal M}_1\cup{\cal M}_3\cup{\cal M}_4\cup{\cal M}_5)x_0)\le\mu({\cal M}_3)$. Thus $\mu({\cal M}_1)=\mu({\cal M}_4)=\mu({\cal M}_5)=0$. 

Part b) implies that $\mu({\cal M}_2)+\mu({\cal M}_7)=\mu({\cal M}_2\cup{\cal M}_7)=\mu(({\cal M}_2\cup{\cal M}_7)x_1)\le\mu({\cal M}_7)$. So $\mu({\cal M}_2)=0$.

Now part c) says that $\mu({\cal M}_7)=\mu({\cal M}_7x_0^{-1})\le\mu({\cal M}_2)$ so $\mu({\cal M}_7)=0$. Finally, it follows from d) that $\mu({\cal M}_3)=\mu({\cal M}_3x_1^{-1})\le\mu({\cal M}_4)$ so $\mu({\cal M}_3)=0$. 

Thus we proved that $\mu({\cal M}_1)=\mu({\cal M}_2)=\mu({\cal M}_3)=\mu({\cal M}_4)=\mu({\cal M}_5)=\mu({\cal M}_7)=0$, which completes the proof.
\vspace{1ex}

An immediate corollary can be extracted from that. Let us say that a subset in $F$ is {\em 1-measured} whenever its complement is zero-measured. It is obvious that a finite union of zero-measured sets is zero-measured. So a finite intersection of 1-measured sets is 1-measured. The property of a set to be 1-measured is invariant under right multiplication. Therefore, for any 1-measured subset ${\cal N}$ and for any positive integer $k$, the intersection $\bigcap\limits_{|w|\le k}{\cal N}w$ will be 1-measured, where the intersection is taken over all words of length $\le k$.

The set ${\cal M}_6$ can be characterized as follows. It consists of those elements whose canonical diagrams are right divisible by a diagram $X_0^{-1}X_1^{-1}$ of two cells. Taking ${\cal N}={\cal M}_6$ and choosing $k\gg1$ large enough, we see that the set of elements whose canonical diagrams are right divisible by any fixed negative diagram right divisible by $X_0^{-1}X_1^{-1}$, form a 1-measured subset. According to Proposition~\ref{remzero}, one can assume without loss of generality that Folner sets consist of the above diagrams (if $F$ is amenable). This can help in clarifying the structure of Folner sets whenever they exist.

\end{document}